\documentclass[11pt]{article}
\usepackage{amsmath}

\newcommand{\ignore}[1]{}

\newcounter{subscript}

\begin{document}
\begin{center}
{\bf \large Using a computer algebra system to simplify}
\end{center} 
\vspace{0.03in}
\begin{center}
{\bf \large expressions for Titchmarsh-Weyl m-functions}
\end{center} 
\vspace{0.03in}
\begin{center}
{\bf \large associated with the Hydrogen Atom on the half line}
\end{center} 

\vspace{1.0in}
\begin{center}
by
\end{center}
\vspace{1.0in}
\begin{center}
\LARGE Cecilia \LARGE Knoll  and \LARGE Charles \LARGE Fulton
\end{center}

\vspace{2.0in}
\begin{center}
{\huge \bf TECHNICAL REPORT} \\
\vspace{0.3in}
{\LARGE \bf  FLORIDA INSTITUTE  OF  TECHNOLOGY} \\
\vspace{0.1in}
{\LARGE \bf  September 2007}
\end{center}

\pagebreak

\title{Using a computer algebra system to simplify expressions for Titchmarsh-Weyl
m-functions associated with the Hydrogen Atom on the half line\footnote{This research partially supported by National Science 
Foundation Grant DMS-0109022 to Florida Institute of Technology.}}

\date{  }

\maketitle

\noindent
\author{CECILIA KNOLL \\ Department of Applied Mathematics \\
Florida Institute of Technology \\ Melbourne, Florida 32901-6975 \\ \\
CHARLES FULTON \\ Department of Mathematical Sciences \\
Florida Institute of Technology \\ Melbourne, Florida 32901-6975 \\ \\
}

\begin{abstract}
{\bf Abstract:} In this paper we give simplified formulas for certain polynomials
which arise in some new Titchmarsh-Weyl m-functions for the radial part of the
separated Hydrogen atom on the half line $(0,\infty)$ and two independent programs
for generating them using the symbolic manipulator Mathematica.
\end{abstract}

\section{Introduction}                                     
\setcounter{equation}{0}

Recently Fulton\cite{FULTON} and Fulton and Langer \cite{FL} considered the following Sturm-Liouville problem for the hydrogen atom on the half line, $x \in (0, \infty)$:
\begin{equation}\label{1.1}
 -y''+\left( -\frac{a}{x} +\frac{\ell
(\ell+1)}{x^2}\right) y = \lambda y , \,\,\,\, 0<x<\infty
\end{equation}
\begin{equation}
\label{1.2} \lim_{x\rightarrow 0} W_x \left( y, x^{\frac{1}{2}} \,\,
J_{2\ell+1} (\sqrt{4ax}) \right)=0.
\end{equation}

In \cite{FULTON} a fundamental system of Frobenius solutions defined at $x=0$ was introduced having the forms:
\begin{equation}\label{1.3}
\phi(x, \lambda)=x^{\ell+1}\left[1+\sum^\infty_{n=1}a_n(\lambda)x^n\right]=\frac{1}{(-2\sqrt{\lambda})^{\ell+1}}M_{\beta, \ell+\frac{1}{2}}(-2ix\sqrt{\lambda}), \quad \beta :=\frac{ia}{2\sqrt{\lambda}}
\end{equation}

\pagebreak
and
\begin{equation}\label{1.4}
\theta(x, \lambda)=-\frac{1}{2\ell+1}\left[K_\ell(\lambda)\phi(x, \lambda)\ln x + x^{-1}+\sum^\infty_{n=1}d_n(\lambda)x^{-\ell+n}\right]
\end{equation}
Where $a_n(\lambda)$ and $d_n(\lambda)$ are polynomials in $\lambda$, $M_{\beta, \ell+\frac{1}{2}}$ is the Whittaker function of first kind and
\begin{equation}\label{1.5}
K_\ell(\lambda)=-\frac{a}{(2\ell+1)!(2\ell)!}\prod^\ell_{j=1}\left[4\lambda j^2 + a^2\right] = \frac{(-2i\sqrt{\lambda})^{2\ell+1}(-\ell-\beta)_{2\ell+1}}{(2\ell)!(2\ell+1)!}.
\end{equation}
Then a Titchmarsh-Weyl $m$-function was introduced in \cite{FULTON} by the requirement
\begin{equation}\label{1.6}
\theta(x, \lambda)-m_\ell(\lambda)\phi(x, \lambda) \in L_2(0, \infty),
\end{equation}
which gives $m_\ell$ as
\begin{eqnarray}\label{1.7}
m_\ell(\lambda)&=&-a k_\ell(\lambda)\left\{\log(-2i\sqrt{\lambda})+\Psi\left(1-\frac{ia}{2\sqrt{\lambda}}\right)-H_{2\ell}+2\gamma\right\}\nonumber\\
&-& ak_\ell(\lambda)\left\{\sum^\ell_{j=1}\frac{1}{\left(j-\frac{ia}{2\sqrt{\lambda}}\right)}\right\} + \frac{(i\sqrt{\lambda})^{2\ell+1}}{(2\ell+1)!}\left[\sum^{2\ell}_{k=0}\frac{2^k\left(-\ell-\frac{ia}{2\sqrt{\lambda}}\right)_k}{k!(2\ell+1-k)}\right]
\end{eqnarray}
where $\Psi (z) =\frac{\Gamma'(z)}{\Gamma(z)}$ is the psi or digamma function, $H_{2\ell}=\sum^{2\ell}_{j=1}\frac{1}{j}$, $\gamma =$ Euler's constant, and the Pochhammer symbol is defined for any complex $z$ as $(z)_k=z(z+1)...(z+k-1)$. Here $k_\ell(\lambda)$ is the polynomial of degree $\ell$ defined by 
\begin{equation}\label{1.8}
k_\ell(\lambda)=-\frac{1}{a(2\ell+1)} K_\ell(\lambda) = -\frac{1}{(2\ell+1)!}\prod^{\ell}_{j=1} \left[4\lambda j^2 + a^2\right].
\end{equation}
For $\ell=0$, we define $k_0(\lambda)=1$. The function $m_\ell(\lambda)$ is an analytic function in the half planes Im $\lambda<0$ and Im $\lambda >0$, has poles on the negative $x$-axis at the eigenvalues of the problem \eqref{1.1}-\eqref{1.2}, and a branch cut, corresponding to the continuous spectrum, on the positive real $\lambda$-axis. In \cite{FL} a Pick-Nevalinna representation of $m_0(\lambda)$ was obtained for $\ell=0$ and for $\ell \ge 1$ it was shown that $m_\ell(\lambda)$ has a Q-function representation which puts it in the class $N_\kappa$ of generalized Nevalinna functions with $\kappa=\left[\frac{\ell+1}{2}\right]$.

Our purpose in this paper is to show that the last two terms in $m_\ell(\lambda)$ can be decomposed into real and imaginary parts as 
\begin{equation}\label{1.9}
-a k_\ell(\lambda)\left\{\sum^{\ell}_{j=1} \frac{1}{j-\frac{ia}{2\sqrt{\lambda}}}\right\} + \frac{(i\sqrt{\lambda})^{2\ell+1}}{(2\ell+1)!}\left[\sum^{2\ell}_{k=0}\frac{2^k\left(-\ell-\frac{ia}{2\sqrt{\lambda}}\right)_k}{k!(2\ell+1-k)}\right]=i\sqrt{\lambda} k_\ell(\lambda)+\frac{r_\ell(\lambda)}{2\ell+1},
\end{equation}
where $r_\ell(\lambda)$ is a polynomial of degree $\ell$ in $\lambda$.

Effectively, the decomposition \eqref{1.9} into real and imaginary parts becomes a defining equation for $r_\ell(\lambda)$. Our first {\it Mathematica} program enables the decomposition to be verified. Next, in Section 3 we separate the left hand side of \eqref{1.9} into real and imaginary parts by introducing polynomial representations of $(-\ell-t)_k$ in $t$ and 
$$
\underset{j \ne m}{\prod^{\ell}_{j=1}}\left(\lambda+\frac{a^2}{4j^2}\right) \text{ in }\lambda,
$$
so that the real part can be represented as a polynomial of degree $\ell$. This yields a real, somewhat explicit, representation for $r_\ell/(2\ell+1)$, and using it a second {\it Mathematica} program shows that this real representation yields the same result as the first {\it Mathematica} program.

\section{The First {\it Mathematica} program}
\setcounter{equation}{0}
A simple program in {\it Mathematica} can be used to verify that the polynomial $r_\ell(\lambda)$ in \eqref{1.9} is real valued. This makes use of the built-in function {\bf Pochhammer [$a,k$]} which executes the multiplications in the Pochhammer symbol $(a)_k$.

Using the second expression in \eqref{1.5} for $K_\ell(\lambda)$, we have using \eqref{1.8} that 
$$
b:=-a k_\ell(\lambda)= \frac{K_\ell(\lambda)}{2\ell+1}=\frac{(-2i\sqrt{\lambda})^{2\ell+1}(-\ell-\beta)_{2\ell+1}}{[(2\ell+1)!]^2}.
$$
Accordingly, solving \eqref{1.9} for $\frac{r_\ell(\lambda)}{2\ell+1}$, we have 
\begin{equation}\label{2.1}
\frac{r_\ell(\lambda)}{2\ell+1}=b\left\{\sum^\ell_{j=1}\frac{1}{j-\beta}\right\} + c + \frac{bi\sqrt{\lambda}}{a}, \text{ where } c:=\frac{(i\sqrt{\lambda})^{2\ell+1}}{(2\ell+1)!}\left[\sum^{2\ell}_{k=0}\frac{2^k\left(-\ell-\frac{ia}{2\sqrt{\lambda}}\right)_k}{k!(2\ell+1-k)}\right].
\end{equation}

Following is the {\it Mathematica} program which implements equation \eqref{2.1}. The output for $\frac{r_\ell(\lambda)}{2\ell+1}$ for $\ell=1,2,3,$ and 4 is shown. The program was executed up to $\ell=30$ showing that $r_\ell(\lambda)$ remained real-valued. Observe that the constant $a$ multiplies all terms of $\frac{r_\ell(\lambda)}{2\ell+1}$ which in the {\it Mathematica} output is ans1. This is also poved in \eqref{3.7} below.

\begin{center}{\bf First Program}\end{center}
\noindent{\bf Program input: $\ell$}
\begin{equation}
\begin{array}{l}
\ell=4;\\
\text{While}[\ell<6, \beta =\frac{Ia}{2\sqrt{\lambda}};\\
b=\text{Expand}[(-2I\sqrt{\lambda})^{2l+1}\frac{\text{Pochhammer}[-\ell-\beta, 2\ell+1]}{(2\ell+1)!^2}; \text{(*This generates $K_\ell/(2\ell+1)$. See (1.5).*)}\\
c=\text{Apart}[\text{Simplify}[\frac{(I\sqrt{\lambda})^{2\ell+1}}{(2\ell+1)!}\sum^{2\ell}_{k=0}\frac{\text{Pochhammer}[-\ell-\beta, k]2^k}{k!(2\ell+1-k)}]];\\
g=b\left(-\sum^{2\ell}_{j=1}\frac{1}{j}\right);\\
e=b\sum^{\ell}_{j=1}\frac{1}{j-\beta};\\
f=\frac{bI\sqrt{\lambda}}{a};\\
\text{Print[``For $\ell= $", $\ell$]};\\
ans1=\text{Simplify}[c+e+f];\\
\text{Print}[``ans1=" ans1]; \text{(*This generates the RHS of (2.1).*)}\\
ans2=\text{Series}[ans1, \{\lambda,0,\ell\}];
\text{Print}[``=", ans2];\\
ans3=\text{Simplify}[g+c+e+f];\text{Print}[``ans3=", ans3]; \text{(*This generates the RHS of \eqref{2.2}*)}\\
ans4=\text{Series}[ans3, \{\lambda, 0, \ell\}];\text{Print}[``=", ans4];\\
\ell=\ell+1]
\end{array}\nonumber
\end{equation}

\noindent{\bf Program output:}\\
For $\ell=1$
\begin{eqnarray}
ans1&=&-\frac{a\lambda}{36}\nonumber\\
ans3&=&\frac{1}{72}(3a^3+10a\lambda)\nonumber = \frac{a^3}{24}+\frac{5a\lambda}{36}\nonumber\\
\end{eqnarray}

For $\ell=2$
\begin{eqnarray}
ans1&=&-\frac{a\lambda(a^2+13\lambda)}{7200} = -\frac{a^3\lambda}{7200}-\frac{13a\lambda^2}{7200}\nonumber\\
ans3&=&\frac{25a^5+476a^3\lambda+1288a\lambda^2}{172800} = \frac{a^5}{6912}+\frac{119a^3\lambda}{43200}+\frac{161a\lambda^2}{21600}\nonumber\\
\end{eqnarray}

For $\ell=3$
\begin{eqnarray}
ans1&=&-\frac{a\lambda(a^4+46a^2\lambda+400\lambda^2)}{8467200} = -\frac{a^5\lambda}{8467200}-\frac{23a^3\lambda^2}{4233600}-\frac{a\lambda^3}{21168}\nonumber\\
ans3&=&\frac{49a^7+2684a^5\lambda+35656a^3\lambda^2+88896a\lambda^3}{508032000}\nonumber\\
&=&\frac{a^7}{10368000}+\frac{671a^3\lambda}{127008000}+\frac{4457a^3\lambda^2}{63504000}+\frac{463a\lambda^3}{2646000}\nonumber\\
\end{eqnarray}

For $\ell=4$
\begin{eqnarray}
ans1&=&-\frac{a\lambda(a^6+107a^4\lambda+3124a^2\lambda^2+22548\lambda^3)}{32920473600}\nonumber\\
&=&-\frac{a^7\lambda}{32920473600}-\frac{107a^5\lambda^2}{32920473600}-\frac{781a^3\lambda^3}{8230118400}-\frac{1879a\lambda^4}{2743372800}\nonumber\\
ans3&=&\frac{761a^9+90200a^7\lambda+3204208a^5\lambda^2+36438400a^3\lambda^3+86960256a\lambda^4}{36870930432000}\nonumber\\
&=&\frac{761a^9}{36870930432000}+\frac{451a^7\lambda}{184354652160}+\frac{4087a^3\lambda^3}{11522165760}+\frac{226459a\lambda^4}{96018048000}.\nonumber\\
\end{eqnarray}
Output for the polynomials
\begin{equation}\label{2.2}
p_\ell(\lambda):=\frac{r_\ell(\lambda)}{2\ell+1}+a k_\ell(\lambda)H_{2\ell},
\end{equation}
which arise in the representation of the $m_\ell$ function from \cite[Equation(8.15)]{FULTON},
$$
m_\ell(\lambda)=k_\ell(\lambda)\left\{-a\log(-2i\sqrt{\lambda})-a\Psi\left(1-\frac{ia}{2\sqrt{\lambda}}\right)-2\gamma a +i\sqrt{\lambda}\right\}+p_\ell(\lambda)
$$
is also given below for $\ell=1,2,3,4.$

\begin{eqnarray}
p_1 &=& \frac{a^3}{24} + \frac{5a\lambda}{36}, \nonumber \\
p_2 &=& \frac{a^5}{6912}
+ \frac{119 a^3}{43200} 
\lambda + \frac{161a}{21600}
\lambda^2, \nonumber\\
p_3 & = &  
\frac{a^7}{10368000} + 
\frac{671 a^5}{127008000}
\lambda +  \frac{4457
a^3}{63504000} \lambda^2 +  
\frac{463 a}{2646000} \lambda
^3\nonumber\\
p_4& = &  \frac{761
a^9}{36870930432000} +
\frac{451a^7}{184354652160} 
\lambda +  \frac{4087 a^5}{47029248000} \lambda^2 + \frac{11387
a^3}{11522165760} 
\lambda^3+ 
\frac{226459a}{96018048000}
\lambda^4 .\nonumber
\end{eqnarray}
Here $p_\ell$ is printed as $ans3$ in the above program.

\section{Explicit Representation for $r_{\ell}$}
\setcounter{equation}{0}
In this section we give a method for separating the expression in equation \eqref{1.9} into real and imaginary parts, yielding a real representation for the polynomial $r_\ell(\lambda)$. The difficulty arises from the complicated product in the Pochhammer symbol $(\ell-\beta)_k$ where $\beta=\frac{ia}{2\sqrt{\lambda}}$. Replacing $\beta$ by a real variable $t$, we let the coefficients of the polynomial $(-\ell-t)_k$ be defined by
\begin{equation}\label{3.1}
g_k(t):=(-\ell-t)_k = \prod^{k-1}_{j=0}(-\ell+j-t)=\sum^k_{n=0}\alpha(k,n)t^n=\sum^{k_1}_{j=0}\alpha(k, 2j)t^{2j}+\sum^{k_2}_{j=0}\alpha(k, 2j+1)t^{2j+1},
\end{equation}
where $k_1=\left[\frac{k}{2}\right]$ and $k_2=\left[\frac{k-1}{2}\right]$.
Here $\alpha(k,n)=\alpha_\ell(k, n)$, and we are interested for fixed $\ell$ to have $\alpha_\ell(k, n)$ available for all $0\le k \le 2\ell$, and all $0 \le n \le k$. For $n=0$, the constant term is $\alpha_\ell(k, 0)=\prod^{k-1}_{j=0}(-\ell+j)=(-\ell)_k$, $k = 0,1, \cdots 2\ell$. Formulas for larger values of $n$ become increasingly more complicated and are not known in closed form. We can, however, represent the real and the imaginary parts of \eqref{1.9} in terms of $\alpha(k,n)$. Putting $t=\beta=\frac{ia}{2\sqrt{\lambda}}$ in \eqref{3.1} gives
\begin{eqnarray}\label{3.2}
g_k(\beta)&=&\sum^k_{n=0}\alpha(k, n)\left(\frac{a}{2}\right)^n i^n (\lambda^{-\frac{1}{2}})^n\nonumber\\
&=&\sum^{k_1}_{j=0}(-1)^j \alpha(k, 2j)\left(\frac{a}{2}\right)^{2j}i^n(\lambda^{-\frac{1}{2}})^{2j}\nonumber\\
&&+i\sum^{k_2}_{j=0}(-1)^j\alpha(k, 2j+1)\left(\frac{a}{2}\right)^{2j+1}(\lambda^{-\frac{1}{2}})^{2j+1}.
\end{eqnarray}

Now for the $(i\sqrt{\lambda})^{2\ell+1}$ term in \eqref{1.9} we have $(i)^{2\ell+1}=(-1)^\ell i$, $\ell=0,1,2,...$ Thus the second sum in \eqref{1.9} may be written as
\begin{eqnarray}\label{3.3}
&& i(-1)^\ell \sum^{2\ell}_{k=0}\frac{2^k}{(2\ell+1)!k!(2\ell+1-k)} g_k(\beta)(\lambda^{\frac{1}{2}})^{2\ell+1}\nonumber\\
&=&-(-1)^\ell \sum^{2\ell}_{k=1} \frac{2^k}{(2\ell+1)!k!(2\ell+1-k)}\left[\sum^{k_2}_{j=0}(-1)^j\alpha(k, 2j+1)\left(\frac{a}{2}\right)^{2j+1}(\lambda^{\frac{1}{2}})^{2\ell+1-(2j+1)}\right]\nonumber\\
&&+ i(-1)^\ell  \sum^{2\ell}_{k=0}\frac{2^k}{(2\ell+1)!k!(2\ell+1-k)}\left[\sum^{k_1}_{j=0}(-1)^j \alpha(k, 2j)\left(\frac{a}{2}\right)^{2j+1}(\lambda^{\frac{1}{2}})^{2\ell+1-2j}\right]\nonumber\\
&=& \sum^{2\ell}_{k=0}\frac{2^k A_k}{(2\ell+1)!k!(2\ell+1-k)}+i\sum^{2\ell}_{k=0}\frac{2^k B_k}{(2\ell+1)!k!(2\ell+1-k)}
\end{eqnarray}
where 
\begin{equation}\label{3.4}
A_k=(-1)^{\ell+1}\sum^{k_2}_{j=0}(-1)^j\alpha(k, 2j+1)\left(\frac{a}{2}\right)^{2j}\lambda^{\ell-j} 
\end{equation}
and
\begin{equation}\label{3.5}
B_k = (-1)^\ell\sum^{k_1}_{j=0} (-1)^j \alpha(k, 2j)\left(\frac{a}{2}\right)^{2j}\lambda^{\ell-j}\lambda^{\frac{1}{2}}
\end{equation}
Similarly, we may decompose the first term in \eqref{1.9} as
\begin{equation}\label{3.6}
\sum^\ell_{m=1}\frac{1}{m-\frac{ia}{2\sqrt{\lambda}}}=\sum^{\ell}_{m=1}\frac{4\lambda m}{4 \lambda m^2 + a^2}+ i \sum^{\ell}_{m=1} \frac{2a \lambda^{\frac{1}{2}}}{4\lambda m^2+a^2}.
\end{equation}
Using \eqref{1.8} and cancelling one $4\lambda j^2+a^2$ factor for each term in the above sums, we find that
\begin{eqnarray}\label{3.7}
-ak_\ell(\lambda)\left[\sum^{\ell}_{m=1}\frac{1}{m-\frac{ia}{2\sqrt{\lambda}}}\right]&=&-\frac{a}{[(2\ell+1)!]^2}\left[\sum^{\ell}_{m=1}4\lambda m\underset{j \ne m}{\prod^\ell_{j=1}}(4\lambda j^2+a^2)\right]\nonumber\\
&&+i \left(-\frac{a}{[(2\ell+1)!]^2}\left[\sum^{\ell}_{m=1}2a\lambda^{\frac{1}{2}}\prod^\ell_{j=1}(4\lambda j^2 + a^2)\right]\right).
\end{eqnarray}
Combining \eqref{3.3} and \eqref{3.7} and taking real and imaginary parts of \eqref{1.9}, we thus obtain:
\begin{equation}\label{3.8}
\sqrt{\lambda}k_\ell(\lambda)=-\frac{2a^2\sqrt{\lambda}}{[(2\ell+1)!]^2}\left[\sum^{\ell}_{m=1}\underset{j\ne m}{\prod^{\ell}_{j=1}}(4\lambda j^2 +a^2)\right]+\sum^{2\ell}_{k=0}\frac{2^kB_k}{(2\ell+1)!k!(2\ell+1-k)}
\end{equation}
and
\begin{equation}\label{3.9}
\frac{r_\ell(\lambda)}{2\ell+1}=-\frac{4\lambda a}{[(2\ell+1)!]^2}\left[\sum^{\ell}_{m=1}m\underset{j \ne m}{\prod^{\ell}_{j=1}}(4\lambda j^2 +a^2)\right]+\sum_{k=1}^{2\ell}\frac{2^k A_k}{(2\ell+1)!k!(2\ell+1-k)}.
\end{equation}
Since every $A_k$ involves $a$ as a factor, it follows that a factor of $a$ occurs in both terms on the right hand side of \eqref{3.9}. Hence $a=0$ will give $r_\ell(\lambda)=0$.

To bring the expression in \eqref{3.9} into a simpler form, we now isolate the coefficients of the powers of $\lambda$. For the second term in \eqref{3.9} we insert the formula \eqref{3.4} for $A_k$ into the sum, change the summation index $(m=\ell-j)$, and then interchange the order of summation:
\begin{eqnarray}\label{3.10}
&&\sum^{2\ell}_{k=1}\frac{2^k A_k}{(2\ell+1)!k!(2\ell+1-k)}\nonumber\\
&=&\sum^{2\ell}_{k=1}\quad\sum^{k_2}_{j=0}\frac{2^k(-1)^{\ell+1+j}\alpha(k, 2j+1)\left(\frac{a}{2}\right)^{2j+1}\lambda^{\ell-j}}{(2\ell+1)!k!(2\ell+1-k)}\nonumber\\
&=&\sum_{k=1}^{2\ell}\quad\sum^{\ell}_{m=\ell-k_2}\left[\frac{2^k(-1)^{m+1}\alpha(k, 2(\ell-m)+1)\left(\frac{a}{2}\right)^{2(\ell-m)+1}}{(2\ell+1)!k!(2\ell+1-k)}\right]\lambda^m\nonumber\\
&=&\sum^{\ell}_{m=1}\quad\sum^{2\ell}_{k=2(\ell-m)+1}\left[\frac{2^k(-1)^{m+1}\alpha(k, 2(\ell-m)+1)\left(\frac{a}{2}\right)^{2(\ell-m)+1}}{(2\ell+1)!k!(2\ell+1-k)}\right]\lambda^m\nonumber\\
&=&\sum^\ell_{j=1}d_j\lambda^j
\end{eqnarray}
where

\begin{equation}\label{3.11}
d_j:=\frac{a}{2}\sum^{2\ell}_{k=2(\ell-j)+1}\left[\frac{2^k(-1)^{m+1}\left(\frac{a}{2}\right)^{2(\ell-j) }\alpha(k, 2(\ell-j)+1)}{(2\ell+1)!k!(2\ell+1-k)}\right].
\end{equation}
Similarly, to isolate the powers of $\lambda$ in the first term in \eqref{3.7} we first define the coefficients $\gamma(m,n)$ of the $\ell-1$ degree polynomial,
\begin{equation}\label{3.12}
\underset{j \ne m}{\prod^{\ell}_{j=1}}\left(\lambda+\frac{a^2}{4j^2}\right)=\sum^{\ell-1}_{n=0}\gamma(m,n)\lambda^n.
\end{equation}
Insertion of this into the first term in \eqref{3.9}, changing the summation index $(j=n+1)$, and interchanging the order of summation then yields:
\begin{eqnarray}\label{3.13}
&&\frac{-4\lambda a}{[(2\ell+1)!]^2}\left[\sum^{\ell}_{m=1} m \underset{j \ne m}{\prod^{\ell}_{j=1}} (4\lambda j^2 + a^2)\right]\nonumber\\
&=&\frac{-a4^\ell(\ell!)^2}{[(2\ell+1)!]^2}\left[\sum^{\ell}_{m=1}\sum^{\ell-1}_{n=0}\left[\frac{\gamma(m,n)}{m}\right]\lambda^{n+1}\right]\nonumber\\
&=&\frac{-a4^\ell(\ell!)^2}{[(2\ell+1)!]^2}\left[\sum^{\ell}_{m=1}\sum^{\ell}_{j=1}\left[\frac{\gamma(m,j-1)}{m}\right]\lambda^{j}\right]\nonumber\\
&=&\frac{-a4^\ell(\ell!)^2}{[(2\ell+1)!]^2}\left[\sum^{\ell}_{j=1}\sum^{\ell}_{m=1}\left[\frac{\gamma(m,j-1)}{m}\right]\lambda^{j}\right]\nonumber\\
&=&\sum^{\ell}_{j=1}c_j\lambda^j
\end{eqnarray}
where
\begin{equation}\label{3.14}
c_j := -\frac{a4^\ell(\ell!)^2}{[(2\ell+1)!]^2} \sum^{\ell}_{m=1}\left(\frac{\gamma(m,j-1)}{m}\right).
\end{equation}
Putting \eqref{3.10} and \eqref{3.13} in \eqref{3.9} we thus have for the polynomial $r_\ell(\lambda)/(2\ell+1)$ the representation,
\begin{equation}\label{3.15}
\frac{r_\ell(\lambda)}{2\ell+1}=\sum^\ell_{j=1}(c_j+d_j)\lambda^j.
\end{equation}
Here it is clear that $a$ is a common factor in $c_j$ and $d_j$ for all $j$, and hence for the case $a=0$, $r_\ell\equiv 0$. The expression \eqref{3.15} is somewhat more explicit than the expression \eqref{2.1} and \eqref{3.9} since the coefficients of the powers of $\lambda$ are isolated, and there are no complex terms present. On the other hand, further simplifications are certainly desirable; however, this requires closed form formulas for $\alpha(k,n)$ and $\gamma(m,n)$ which remain elusive. However, the formulas for $c_j$ and $d_j$ are easily implemented using a symbolic manipulator.

\section{The Second {\it Mathematica} Program}
\setcounter{equation}{0}
As an independent check on the first {\it Mathematica} program we implemented the formulas \eqref{3.11} and \eqref{3.14} for $d_j$ and $c_j$ and computed the polynomial of degree $\ell$ in \eqref{3.15}. This required computing and storing the coefficient $\alpha(k,n)$ and $\gamma(m,n)$ in \eqref{3.1} and \eqref{3.12}. Following is the {\it Mathematica} program that does this to compute $r_\ell(\lambda)/(2\ell+1)$. The output  for $\ell=1,2,3$ and 4 is shown. The program was executed up to $\ell=30$ and gave exact agreement with the first {\it Mathematica} program.

\begin{center}
{\bf Second Program}
\end{center}

\noindent{\bf Program input: $\ell$}
\begin{equation}\begin{array}{l}
\ell=4; k=0;\nonumber\\
\text{While} [k \le 2\ell,\nonumber\\
\quad p[t\_] = \text{Expand}[\text{Pochhammer}[-\ell -t, k]]; n=k; i=0;\nonumber\\
\quad\text{While}[i\le n,\nonumber\\
\quad\quad\alpha[k,i]=p[0];\nonumber\\
\quad\quad p[t\_]=\text{Simplify}[\text{Expand}[(p[t]-\alpha[k,i])/t]];\nonumber\\
\quad\quad i++]\nonumber\\
\quad k++]\nonumber\\
\text{Clear}[k]; j=1;\nonumber\\
\text{While}[j\le\ell,\nonumber\\
\quad d[j]=\sum^{2\ell}_{k=2(\ell-j)+1}\frac{2^k(-1)^{j-1}\alpha[k, 2\ell-2j+1]\left(\frac{a}{2}\right)^{2\ell-2j+1}}{(2\ell+1)!k!(2\ell+1-k)}; j++]\nonumber\\
\text{Clear}[j]; m=1;\nonumber\\
\text{While}[m\le\ell,\nonumber\\
\quad q[t\_]=\left(\prod^\ell_{j=1}\left(t+\frac{a^2}{4j^2}\right)\right)/\left(t+\frac{a^2}{4m^2}\right); n=\ell-1; i=0;\nonumber\\
\quad\text{While}[i\le n,\nonumber\\
\quad\quad\gamma[m,i]=q[0];\nonumber\\
\quad\quad q[t\_]=\text{Simplify}[\text{Expand}[(q[t]-\gamma[m,i])/t]];\nonumber\\
\quad\quad i++]\nonumber\\
\quad m++]\nonumber\\
\text{Clear}[i]; j=1;\nonumber\\
\text{While}[j \le \ell,\nonumber\\
\quad c[j]=\frac{-a4^\ell(\ell!)^2}{((2\ell+1)!)^2}\sum^{\ell}_{m=1}\frac{\gamma[m, j-1]}{m}; j++;\nonumber\\
\text{Print}[\sum^{\ell}_{j=1}(c[j]+d[j])t^j]\nonumber
\end{array}\end{equation}

\noindent{\bf Output of second program}\\
For $\ell=1:$
$$
\frac{r_\ell}{(2\ell+1)}=-\frac{at}{36}
$$
For $\ell=2:$
$$
\frac{r_\ell}{(2\ell+1)}=-\frac{a^3t}{7200}-\frac{13 a t^2}{7200}
$$
For $\ell=3:$
$$
\frac{r_\ell}{(2\ell+1)}=-\frac{a^5 t}{8467200}-\frac{23a^3t^2}{4233600}-\frac{at^3}{21168}
$$
For $\ell=4:$
$$
\frac{r_\ell}{(2\ell+1)}=-\frac{a^7 t}{32920473600}-\frac{107a^5t^2}{32920473600}-\frac{781a^3t^3}{8230118400}-\frac{1879at^4}{2743372800}
$$
This output is in agreement with the output of the first program, and was also checked up to $\ell=30$. The two {\it Mathematica} programs verify that the imaginary part of the expression \eqref{2.1} for $r_\ell(\lambda)/(2\ell+1)$ is zero for all values of $\ell$ for which the programs were executed. A general proof that this is true for all $\ell$, that is, a proof of \eqref{3.8} (or, equivalently, a proof of \eqref{3.9}) requires that the $\alpha(k,n)$ and $\gamma(m,n)$ coefficients be obtained in a simplified form, and this appears to be a formidable task. It is quite difficult to perform induction on $\ell$ to prove \eqref{3.8} or \eqref{3.9} because of the complicated manner in which $\alpha(k,n)=\alpha_\ell(k,n)$ and $\gamma(m,n)=\gamma_\ell(m,n)$ change with $\ell$. Nevertheless, it may be possible to construct rigorous proofs by making use of some combinatorial analysis. For example there are general forms for the coefficients of a polynomial having $k$ known roots $x_i$, $i=1,2,...,k$. 
Sen and Krishnamurthy \cite{KRISH}, for example, show that 
\begin{equation}\label{4.1}
q_k(x):=\prod^{k}_{j=1}(x-x_j)=x^k+\sum^{k}_{m=1}a_m x^{k-m}
\end{equation}
with
\begin{equation}\label{4.2}
a_m := (-1)^m S_m := (-1)^m \sum x_1x_2...x_m,
\end{equation}
where the sum is taken over all the products of $x_i$ taken $m$ at a time. Thus, for the polynomial $g_k(t)$ of \eqref{3.1} we have
\begin{equation}\label{4.3}
g_k(t)=(-1)^k\prod^{k-1}_{j=0} [t-(j-\ell)]=(-1)^k\prod^{k}_{j=1}[t-(j-1-\ell)]=(-1)^k\left(t^k+\sum^{k}_{m=1}a_m t^{k-m}\right),
\end{equation}
with
\begin{equation}\label{4.4}
a_m=a_\ell (k, m)=(-1)^m\sum t_1 t_2 ... t_m
\end{equation}
where the sum is taken over all the products of $t_j := j-1-\ell$ taken $m$ at a time. But, a general solution for $a_m$ as a function of $k$ and $\ell$, valid for all $\ell$ and all $0 \le k \le 2\ell$, remains elusive. Another idea would be to establish \eqref{3.8} or \eqref{3.9} by induction on $\ell$, but it unfortunately appears difficult to employ the induction hypothesis. So, an analytic proof of \eqref{3.8} or \eqref{3.9} valid for all $\ell \ge 1$ remains an open problem.  

\section{Conclusion}
The simplification of $a_\ell(k, n)$ and a rigorous proof of \eqref{3.8} and/or \eqref{3.9} remain as open problems. We have, however, given in this note two {\it Mathematica} programs which implicitly establish both of these results.


\vskip 18pt
\begin{thebibliography}{99}

\bibitem{FULTON} C. Fulton, Titchmarsh-Weyl m-functions for Second-order Sturm Liouville Problems with two singular
 endpoints, Math. Nachr. 281 (10) (2008), 1418-1475.

\bibitem{FL} C. Fulton and H. Langer, Sturm-Liouville operators with singularities and Generalized Nevanlinna functions,
 Complex Anal. and Oper. Theory, to appear.

\bibitem{KRISH} E.V. Krishnamurthy and S.K. Sen, Numerical Algorithms: Computations in Science and Engineering, Affiliated East West Press,
  New Delhi, 2001.

\end{thebibliography}
\end{document}